\documentclass[12pt,thmsa]{amsart}
\usepackage{amssymb, euscript}

\def\Cal{\mathcal}

\def\I{{\Cal I}}

\def\bbr{{\mathbb R}}

\def\bbc{{\mathbb C}}

\def\diag{{\hbox{\rm diag}}}

\def\cos{{\hbox{\rm cos}}}

\def\vol{{\hbox{\rm vol}}}
\def\min{{\hbox{\rm min}}}

\def\part{\partial}
\def\intl{\int\limits}
\def\b{\beta}

\def\Gam{\Gamma}
\def\Om{\Omega}
\def\a{\alpha}

\def\om{\omega}

\def\del{\delta}
\def\vp{\varphi}

\def\gam{\gamma}
\def\Lam{\Lambda}
\def\sig{\sigma}
\def\lam{\lambda}

\def\e{\varepsilon}

\font\frak=eufm10

\def\fr#1{\hbox{\frak #1}}

\def\frM{\fr{M}}

\def\f0{f_0}
\def\Fc0{\varphi_0}
\def\rn{\bbr^n}

\def\vni{V_{n,i}}

\def\cos{{\hbox{\rm cos}}}

\numberwithin{equation}{section}

\newcommand{\be}{\begin{equation}}
\newcommand{\ee}{\end{equation}}

\newcommand{\bea}{\begin{eqnarray}}
\newcommand{\eea}{\end{eqnarray}}
\newcommand{\Bea}{\begin{eqnarray*}}
\newcommand{\Eea}{\end{eqnarray*}}

\def\kl{K_\ell}
\def\rl{\bbr^\ell}
\def\rn{\bbr^n}
\def\rnl{\bbr^{n-\ell}}
\def\lv{{\boldsymbol \lam}}

\newtheorem{theorem}{Theorem}[section]
\newtheorem{lemma}[theorem]{Lemma}
\newtheorem{definition}[theorem]{Definition}

\newtheorem{corollary}[theorem]{Corollary}
\newtheorem{proposition}[theorem]{Proposition}

\theoremstyle{remark}
\newtheorem{remark}[theorem]{Remark}

\numberwithin{equation}{section}

\begin{document}

\title[ The Lower Dimensional Busemann-Petty Problem]
{The Lower Dimensional Busemann-Petty Problem for Bodies with the
Generalized Axial Symmetry}

\author{Boris Rubin}
\address{
Department of Mathematics, Louisiana State University, Baton Rouge,
LA, 70803 USA}

\email{borisr@math.lsu.edu}

\thanks{The  research was supported in part by the  NSF grant DMS-0556157
and the Louisiana EPSCoR program, sponsored  by NSF and the Board of
Regents Support Fund.}





\begin{abstract}
The lower dimensional Busemann-Petty problem asks, whether
$n$-dimensional
 centrally symmetric convex bodies with smaller $i$-dimensional
 central
 sections necessarily have smaller volumes. For
 $i=1$, the affirmative answer is obvious. If $i>3$, the
 answer is  negative. For
 $i=2$ and $i=3$,  the problem is still
 open, however,  when the body with smaller sections is
 a body of revolution,  the answer  is
 affirmative. The paper contains a complete solution to the problem
  in the more general situation, when the body with smaller sections is invariant
  under rotations, preserving mutually orthogonal
  subspaces of dimensions $\ell$ and $n-\ell$, respectively,
  so that  $i+\ell\le n$. The answer essentially depends on
   $\ell$. The argument relies on the notion of canonical angles
   between subspaces, spherical Radon transforms, properties  of intersection bodies, and the generalized cosine transforms.
\end{abstract}

\maketitle

\section{Introduction}

\setcounter{equation}{0}

Let $G_{n,i}$  be the Grassmann manifold of $i$-dimensional linear
subspaces  of $\bbr^n$, and let $\vol_i ( \cdot)$ denote the
$i$-dimensional volume function.

\vskip 0.2 truecm

\noindent{\bf Question:} Suppose that $i$ is fixed, and let $A$ and
$B$ be arbitrary origin-symmetric (o.s.)
 convex bodies  in $\bbr^n$ satisfying \be\label{cons} \vol_i(A \cap \xi)
\le \vol_i(B \cap \xi) \quad \forall \xi \in G_{n,i}. \ee Does it
follow that \be\label{cons1} \vol_n(A) \le \vol_n(B) \quad \text{\rm
?} \ee

\vskip 0.2 truecm

This question generalizes  the celebrated  Busemann-Petty problem,
corresponding to $i=n-1$ \cite{BP}. The latter has a long history,
and the  answer is really striking. It is ``Yes''  if and only if $n
\le 4$; see \cite{G}, \cite{K3}, \cite{R2}, \cite{Z2}, and
references therein. For  $1\le i\le n-2$, the problem is more
intriguing. We call it  {\it the lower dimensional Busemann-Petty
problem} (LDBP). If
 $i=1$, the implication (\ref{cons}) $\to$ (\ref{cons1}) is obvious
 for all o.s. star bodies without any convexity assumption. In the case $i=2, \; n=4$,  an affirmative
 answer follows from the solution of the usual  Busemann-Petty problem. For
$3<i \le n-1$, a
 negative answer was first given by Bourgain and  Zhang \cite{BZ}.  The
 proof   in \cite{BZ} was corrected in  \cite{RZ}. An alternative
  proof was given in \cite{K4}. In the cases  $i=2$ and $i=3$ for $n>4$, the answer is
  generally unknown, however,  if the body with smaller sections is
 a body of revolution, the answer is
 affirmative; see  \cite{GZ}, \cite{Z1}, \cite{RZ}.
  It is also known
\cite{BZ}, that when $i=2$ and $B$ is a Euclidean ball, the answer
is affirmative provided that $A$ is convex and sufficiently close to
$B$. On the other hand \cite {Mi2}, for
 $i=2$ and $i=3$, there is a small perturbation $A$
of a Euclidean ball,  so that   the implication (\ref{cons}) $\to$
(\ref{cons1}) is true for arbitrary o.s. star body $B$.
Modifications of  the
  Busemann-Petty problem  were studied in \cite{K3},  \cite{RZ},
  \cite{Y},   \cite {Zv}; see also  \cite{G} , where one can find further references.

\vskip 0.3 truecm

\noindent{\bf Main results.} In the present paper we give a complete
solution to the problem stated above, when the body with smaller
sections is invariant
  under orthogonal transformations preserving mutually orthogonal
  subspaces, say, $p$ and $p^\perp$, of  dimensions $\ell$ and $n-\ell$
  satisfying
   $i+\ell\le n, \;  1\le \ell<n$. Let us choose the
coordinate system in $\rn$ so that $
p=\rl=\overset{n}{\underset{j=n-\ell +1}{\oplus}}\,\bbr e_j$ and
$p^\perp =\rnl=\overset{n-\ell}{\underset{j=1}{\oplus}}\,\bbr e_j$,
where
 $ \, e_1, e_2, \ldots , e_n$ are the relevant coordinate unit
 vectors. Without loss of generality, we  assume $\ell \le
 n-\ell$, i.e., $\ell \le n/2$ (otherwise, the coordinate subspaces can
 be renamed). The case  $\ell=1$ corresponds to
   bodies of revolution.

Consider the subgroup  of orthogonal transformations \be\label{kael}
\kl\!=\!\left \{\gam \in O(n): \gam\! =\! \left[\begin{array} {cc}
\a & 0
\\ 0 & \b \end{array} \right], \quad \a  \!\in  \!O(n-\ell), \; \b  \!\in \!
O(\ell) \right \}. \ee  A star body $A$ is $\kl$-symmetric if $\gam
A=A$ for all  $\gam \in \kl$. Clearly, every $\kl$ -symmetric body
$A$ is origin-symmetric, that is $A=-A$ . We  set  $$ x=(x', x'')
\in \rn, \quad x' \in \rnl, \quad x''\in \rl. $$ Every
$\kl$-symmetric body in $\rn$ can be obtained, for instance, if we
take a 2-dimensional body, which is symmetric with respect to
coordinate axes in the plane $(e_1, e_n)$, and rotate it  about the
subspaces $\rl$ and $\rnl$. A typical example is the
$\;(q,\ell)$-ball \be\label{qlb} B^n_{q,\ell}=\{x: |x'|^q +|x''|^q
\le 1\}, \qquad q>0.\ee

The basic idea of our approach is the following. We observe, that
the relative position of a subspace $\xi \in G_{n,i}$ with respect
to the coordinate subspace $ \rl$ is detemined by $m=\min(i,\ell)$
canonical angles $\om_1, \ldots, \om_m$; see, e.g., \cite{C}. We
define \be G_{n,i}^\ell \!=\!\{\xi \in G_{n,i}: \om_1= \ldots =
\om_m\} \ee to be  the submanifold of all $\xi \in G_{n,i}$ such
that all canonical angles between $\xi $ and $\rl$ are equal. The
structure of the set $G_{n,i}^\ell$ can be understood as follows.
Let $ \lam_1=\cos^2 \om_1, \dots , \lam_m=\cos^2 \om_m$. These  are
eigenvalues of the positive semi-definite matrix
\be\label{forg}r=\left\{
 \begin{array} {ll} \tau '{\rm P}_{\rl}\tau & \mbox{if $ i\le \ell$ ,}\\
\sig'{\rm P}_{\xi}\sig & \mbox{if $ i>\ell$},\\
  \end{array}
\right. \ee
 where $\tau $ and $\sig$ denote arbitrarily fixed
 orthonormal frames which
 span $\xi$ and $\rl$, respectively;
$\tau '$, $\sig'$, ${\rm P}_{\rl}$, and ${\rm P}_{\xi}$  stand for
the corresponding transposed matrices and orthogonal projections. We
arrange $\lam_1,\dots , \lam_m $ in non-increasing order and regard
$\lv=(\lam_1,\dots , \lam_m)$  as a point of the simplex
\be\label{prosh} \Lam_m =\{\lv: 1\ge\lam_1\ge \lam_2\ge\dots \ge
\lam_m \ge 0\}.\ee The edge $\{\lam_1= \ldots = \lam_m\}$ of this
simplex corresponds  to $ G_{n,i}^\ell$.

 Our main results are the following.
\begin{theorem}\label{teor111}
 Let  $1\le \ell \le n/2$, $i+\ell \le n$, and let $A$ be a $K_\ell$-symmetric star body  in
$\bbr^n$.

\vskip 0.3truecm

\noindent {\rm(a)} If $\;1\le i\le \ell$, then the implication
\be\label{implij} \vol_i(A \cap \xi) \le \vol_i(B \cap \xi) \quad
\forall \xi \in G_{n,i}^\ell \Longrightarrow \vol_n(A) \le
\vol_n(B)\ee   is true for every o.s. star body $B$.

\vskip 0.3truecm

\noindent {\rm(b)} If $\;i= \ell+1$ or $i=\ell+2$, then
(\ref{implij}) holds for every o.s. star body $B$ provided that $A$
is convex.
\end{theorem}
\begin{theorem}\label{teor112}
If $\;i>\ell+2$, and $B= B^n_{4,\ell}=\{x: |x'|^4 +|x''|^4 \le 1\},$
then there is an infinitely smooth $K_\ell$-symmetric convex body
$A$, such that  $\vol_i(A \cap \xi) \le \vol_i(B \cap \xi)$ for all
$\xi \in G_{n,i}$, but $  \vol_n(A) > \vol_n(B)$.
\end{theorem}

Some comments are in order.

$1^0$. It might be surprising, that  to make a positive conclusion
in Theorem \ref{teor111}, we do not need {\it all} $i$-dimensional
central sections, as  suggested in the original problem. It suffices
to consider only sections  having equal canonical angles with
respect to $\rl$. More advantages of our {\it method of canonical
angles} are described in Remark \ref {rm26}.

$2^0$. The condition $i+\ell\le n$ in Theorem
 \ref{teor111} excludes the situation when
 $\dim (\xi \cap \rl) \ge  1$ for all $\xi \in G_{n,i}$; see Remark
 \ref{Lzd}. We actually assume

in (a): $i\le \min (\ell, n-\ell)$;

in (b): $\ell \le (n-1)/2$, if $\;i= \ell+1$, and $\ell \le
(n-2)/2$, if $\;i= \ell+2$.

\noindent Regarding (a), the situation, when inequalities $i+\ell>
n$ and $i\le \ell$ hold simultaneously,
 is impossible, because in this case $\ell>n/2$, that contradicts  our initial convention.
Regarding (b), a simple examination shows that the following cases,
which are admissible when $i+\ell>n$, are not presented in Theorem
\ref{teor111}:

\vskip 0.3truecm

(i) $n=2\ell$, when $ i= \ell+1$;

(ii) $n=2\ell$ and  $n=2\ell +1$, when $ i= \ell+2$;

\vskip 0.3truecm

 The validity of the implication (\ref{cons}) $\to$
(\ref{cons1}) in (i) and (ii) is  an open problem. After several
attempts to attack it, we have got an impression that the
difficulties here have the same nature  as those in the original
LDBP for $i=2$ and $3$.

 $3^0$. Another intriguing open problem  is to
check the following

\noindent {\bf Conjecture.}  {\it In the case (b) of Theorem
\ref{teor111},  i.e., when $\;i= \ell+1$ or $i=\ell+2$, there exist
a {\it non-convex}  $\kl$-symmetric body $A$ and an o.s. star body
$B$ so that $$ \vol_i(A \cap \xi) \le \vol_i(B \cap \xi) \quad
\text{for all $\xi \in G_{n,i}$ (not only for  $\xi \in
G_{n,i}^\ell$)},$$ but} $\vol_n(A)> \vol_n(B)$;
 cf. \cite [Theorem 8.2.4] {G} for $n=3, \; \ell =1$.

The paper is organized as follows. In Section 2 we obtain new
 lower dimensional representations for the spherical Radon transform of
$K_\ell$-invariant functions; see Theorem \ref{teo31} and Corollary
\ref{cor32}.  These results are used in Section 3 to prove Theorem
\ref{teor111}. Theorem \ref{teor112} is proved in Section 4, where
we invoke some facts on intersection bodies and the generalized
cosine transforms. The concept of intersection body was introduced
 by Lutwak \cite {Lu} and extended by Zhang \cite{Z1} and
Koldobsky \cite {K4} to lower dimensional sections. Useful
information about these objects can be found in \cite{K3},
 \cite{Mi1}, \cite{R3}.

\noindent{\bf Acknowledgements.} The author is grateful to
Professors Alexander Koldobsky, Erwin Lutwak,  Deane Yang, and
Gaoyong Zhang
 for useful  discussions.

\vskip 0.3truecm

\noindent{\bf Notation:} We use  the standard
 notation $O(n)$   and $SO(n)$ for the orthogonal group and the
 special orthogonal group of $\bbr^{n}$ endowed with the  invariant
 probability measure. For $1\le i<n$, we denote by $G_{n,i}$  the
Grassmann manifold of $i$-dimensional  subspaces  $\xi$ of $\Bbb
R^n$; $d\xi$ stands for the $O(n)$-invariant probability measure on
$G_{n,i}$;
 $\;S^{n-1}$ is the unit sphere in $\bbr^n$; $\, \sig_{n-1}=
2\pi^{n/2}/\Gam (n/2) $ is the area of $S^{n-1}$;  $ \, e_1, e_2,
\ldots , e_n$ denote the coordinate unit vectors; $\frM_{n,i}$ is
 the space of real matrices having $n$ rows and $i$ columns. For $X
\in \frM_{n,i}$, $X'$ denotes   the transpose of $X$, $I_i$ is the
identity $i \times i$
 matrix;
 $$\vni= \{\tau \in \frM_{n,i}: \tau' \tau =I_i \}=O(n)/
O(n-i)$$ is
 the Stiefel manifold of orthonormal $i$-frames in $\bbr^n$. For
 $\tau \in \vni$, $\{\tau\}$ denotes the $i$-dimensional subspace
 spanned by $\tau$.   All vectors in $\rn$ are interpreted as column-vectors.

\section{The Spherical Radon Transform of $K_\ell$-Invariant Functions}

For  functions $f(\theta)$ on $S^{n-1}$ and $\varphi (\xi)$ on
$G_{n,i}$, we define the spherical Radon transform $(R_i f)(\xi)$
 and its
dual  $(R_i^* \varphi)(\theta)$  by \be\label{rts}
 (R_i f)(\xi) = \intl_{S^{n-1}\cap\xi} f(\theta) \, d_\xi \theta, \qquad
  (R_i^* \varphi)(\theta) = \intl_{\xi \ni \theta}  \varphi (\xi)  \, d_\theta \xi,
\ee where  measures $d_\xi \theta$ and $ d_\theta \xi$ are
normalized so that $R_i 1=\sig_{i-1}$ and $ R_i^* 1 = 1$. The
corresponding duality relation has the form
 \be\label{dual} \frac{1}{\sig_{i-1}} \intl_{G_{n,i}} (R_if)(\xi) \vp
(\xi) d \xi = \frac{1}{\sig_{n-1}} \intl_{S^{n-1}} f(\theta) (R_i^*
\vp) (\theta) d\theta \ee and is applicable whenever either side is
finite for $f$ and $\vp$ replaced by $|f|$ and $|\vp|$,
respectively; see \cite{He}, \cite{R1}.

In this section we obtain explicit ``lower dimensional'' expressions
for $R_i f$ when $f$ is $K_\ell$-invariant. We remind that \be
\rn=\rnl \oplus \rl, \qquad
\rnl=\overset{n-\ell}{\underset{j=1}{\oplus}}\,\bbr e_j, \qquad
\rl=\overset{n}{\underset{j=n-\ell +1}{\oplus}}\,\bbr e_j,\ee $1\le
\ell \le n-1$, and set
 \be\label{eell}
\sig=[e_{n-\ell +1},\ldots , e_n]=\left[\begin{array} {c} 0 \\
I_\ell \end{array} \right].\ee Every $\theta \in S^{n-1}$ is
represented in bi-spherical coordinates
as \be\label{bsph} \theta=\left[\begin{array} {c} u\, \sin \om \\
v\, \cos \,\om \end{array} \right], \qquad u\in S^{n-\ell-1}, \quad
v\in S^{\ell-1}, \quad 0\le \om \le \frac{\pi}{2},\ee so that $
d\theta =\sin^{n-\ell-1} \om \,\cos^{\ell-1}\om \,du dv d\om$; see,
e.g., \cite{VK}. Clearly,
$\cos^{2}\om=\theta'\sig\sig'\theta=\theta'{\rm P}_{\rl}\theta$,
where ${\rm P}_{\rl}$ denotes the orthogonal projection onto $\rl$.
The following statement is an immediate consequence of (\ref{bsph}).
\begin{lemma}\label{teo1} A function $f$ on  $S^{n-1}$ is $\kl$-invariant if
and only if there is a function  $f_0$ on $[0,1]$ such that
$f(\theta)=f_0(t)$, where $t^{1/2}=(\theta'{\rm
P}_{\rl}\theta)^{1/2}$ is the cosine of the angle between the unit
vector $\theta$ and the coordinate subspace $\rl$. Moreover,  \bea
\intl_{S^{n-1}}f(\theta)\, d\theta&=&c\,\intl_0^{\pi/2}
\sin^{n-\ell-1} \om \,\cos^{\ell-1}\om \,f_0(\cos^{2}\om) \,
d\om\nonumber
\\\label{gn1}&=&\frac{c}{2}\,\intl_0^{1}t^{\ell/2-1}
(1-t)^{(n-\ell)/2-1} f_0 (t) \, dt, \qquad c=\sig_{\ell
-1}\sig_{n-\ell -1}.\eea
\end{lemma}
\begin{theorem}\label{teo31} Let $1\le i,\ell \le n-1$;
 $m=\min(i,\ell)$.
Let $\om_1, \ldots, \om_m$ be  canonical angles between the subspace
$\xi \in G_{n,i}$ and the coordinate plane $\rl$,\be\label{lpolc}
\lv=\diag (\lam_1,\dots , \lam_m), \qquad \lam_1=\cos^2 \om_1, \dots
, \lam_m=\cos^2 \om_m.\ee Suppose that
  $f$ is a $K_\ell$-invariant
 function on $S^{n-1}$, so that $f(\theta)=f_0(t), \, t=\cos^2
 \om$, where $\om$ is the angle between $\theta$ and $\rl$. Then the
 Radon transform $R_if$ has the form $
(R_if)(\xi)=F(\lv)$, where \be\label{F1}
F(\lv)\!=\!\frac{\sig_{i-\ell
-1}}{2}\!\intl_{S^{\ell-1}}\!\frac{dv}{(v'\lv v)^{i/2
-1}}\!\intl_0^{v'\lv v} \!t^{\ell/2 -1}(v'\lv v -t)^{(i-\ell)/2 -1}
f_0(t)\, dt\ee if $\;i>\ell$, and \be F(\lv)=\intl_{S^{i-1}}
f_0(v'\lv v)\, dv\ee if $\; i\le\ell$.
\end{theorem}
\begin{proof}  We set $$p_i=\left[\begin{array} {c} I_i \\
0 \end{array} \right] \in V_{n,i}, \qquad
\{p_i\}=\overset{i}{\underset{j=1}{\oplus}}\,\bbr e_j, \qquad
\sig=\left[\begin{array} {c} 0 \\
I_\ell \end{array} \right] \in V_{n,\ell},
$$ and let
$\rho_\xi \in SO(n)$ be a rotation that takes the subspace $\{p_i\}$
to $\xi \in G_{n,i}$. Then (set $\theta=\rho_\xi \eta$)
$$
(R_if)(\xi)=\intl_{S^{n-1}\cap\xi} f_0(\theta' \sig\sig'\theta) \,
d_\xi \theta=\intl_{S^{i-1}} f_0(\eta'
\rho'_\xi\sig\sig'\rho_\xi\eta) \, d\eta,$$ $S^{i-1}$ being the unit
sphere in $\{p_i\}$. Let \be\label{43sh} u\!=\!\rho'_\xi\sig\!=\!\left[\begin{array} {c} u_1 \\
u_2 \end{array} \right]\!\in \!V_{n,\ell}, \quad u_1\!=\!p'_i
u\!=\!p'_i \rho'_\xi\sig \!\in \!\frM_{i,\ell}, \quad u_2 \!\in \!
\frM_{n-i,\ell}.\ee Then $\eta'u=\eta'u_1$, and we have \be\label
{44}(R_if)(\xi)=\intl_{S^{i-1}} f_0(\eta' uu'\eta) \,
d\eta=\intl_{S^{i-1}} f_0(\eta' u_1u'_1\eta) \, d\eta.\ee

Consider the case  $\;\ell <i$ and write $u_1$ in the form (cf.
\cite[p. 589]{M})
$$u_1=\gam p_\ell\, r^{1/2}, \qquad \gam \in SO(i),
\quad  p_\ell=\left[\begin{array} {c} I_\ell \\
0 \end{array} \right] \in V_{i,\ell},$$ where $r$ is a positive
semi-definite $\ell \times \ell$ matrix defined by \be r=u'_1
u_1=u'p_i p'_i u=\sig'\rho_\xi p_i p'_i\rho'_\xi\sig=\sig'{\rm
P}_\xi \sig.\ee Hence,
$$ (R_if)(\xi)=\intl_{S^{i-1}} f_0(\eta' \gam p_\ell \,r p'_\ell
\,\gam'\eta) \, d\eta=\intl_{S^{i-1}} f_0(\zeta' p_\ell  \,r p'_\ell
\,\zeta) \, d\zeta.$$ Since $\ell<i$, then $\{p_\ell\}\subset
\{p_i\}$, and we can write $\zeta$ in bi-spherical coordinates
$$
\zeta=\left[\begin{array} {c} v\, \cos \psi\\
w\, \sin \psi\end{array} \right], \qquad v\in S^{\ell-1}, \quad w\in
S^{i-\ell-1}, \quad 0\le \psi \le \frac{\pi}{2},$$ so that $ d\zeta
=\cos^{\ell-1}\psi \,\sin^{i-\ell-1} \psi \,dv dw d\psi$. This gives
$p'_\ell \,\zeta=v\, \cos \psi$, and therefore, \bea
(R_if)(\xi)&=&\sig_{i-\ell -1}\intl_{S^{\ell-1}} dv
\intl_0^{\pi/2}f_0 (v'rv\, \cos^{2}\psi)\,\cos^{\ell-1}\psi
\,\sin^{i-\ell-1} \psi \, d\psi\nonumber\\&=& \frac{\sig_{i-\ell
-1}}{2}\!\intl_{S^{\ell-1}}\!\frac{dv}{(v'r v)^{i/2
-1}}\!\intl_0^{v'r v} \!t^{\ell/2 -1}(v'r v -t)^{(i-\ell)/2 -1}
f_0(t)\, dt.\eea Finally, we diagonalize $r=\sig'{\rm P}_\xi \sig$
 by setting $r=\gam'\lv\gam$, where $\gam \in O(\ell)$ and
$\lv=\diag (\lam_1,\dots , \lam_\ell)$. Changing variables, we
obtain (\ref{F1}).

 Consider the case $\ell\ge i$. We replace $u_1$ in (\ref{44})
by  $p'_i\rho'_\xi \sig$ from (\ref{43sh}) and let $\tau \in
V_{n,i}$ be an arbitrary orthonormal $i$-frame in $\xi$. We can
always choose $\rho_\xi$ so that $\rho_\xi p_i=\tau$. Then $u_1
u'_1=p'_i\rho'_\xi \sig \sig'\rho_\xi p_i=\tau '\sig \sig'\tau $.
The $i \times i$ matrix $s=\tau '\sig \sig'\tau $ is positive
semi-definite and can be diagonalized as above. Hence, (\ref{44})
yields
$$
(R_if)(\xi)=\intl_{S^{i-1}} f_0(\eta's \eta) \,
d\eta=\intl_{S^{i-1}} f_0(\eta'\lv \eta) \, d\eta,$$ as desired.
\end{proof}
\begin{corollary}\label{cor32} If all
canonical angles in Theorem \ref{teo31} are equal, that is,
$\lam_1=\dots = \lam_m=\lam$, then $(R_if)(\xi)=F(\lam)$, where
\be\label{elli} F(\lam)\!=\!\frac{\sig_{i-\ell -1}\,\sig_{\ell
-1}}{2\lam^{i/2
 -1}}\intl_0^\lam \!t^{\ell/2 -1}(\lam -t)^{(i-\ell)/2 -1} f_0(t)\,
 dt\ee
 if $\;i>\ell$, and \be\label{gega}  F(\lam)=\sig_{i -1}\,f_0(\lam) \ee
if $\;i\le\ell$.
\end{corollary}

\begin{remark}\label {Lzd} If $i+\ell >n$, then every $\xi \in G_{n,i}$ has at
least one-dimensional intersection with $\rl$. It means that  some
canonical angles between $\xi$ and $\rl$ are necessarily zero and
therefore, some of the eigenvalues $\lam_1, \ldots \,  \lam_m$ equal
$1$. It follows that for $i+\ell >n$, equalities (\ref{elli}) and
(\ref{gega}) are available only for $\lam =1$. This situation is not
favorable for our purposes, because we will need (\ref{elli}) and
(\ref{gega}) to be available {\it for all} $\lam \in (0,1)$. The
latter is guaranteed if $i\le n-\ell$, when we have ``sufficiently
many'' $i$-dimensional subspaces with the property $\dim (\xi \cap
\rl)=0$.
\end{remark}

Corollary \ref {cor32}  motivates the following
\begin{definition}\label{defi}  We denote by $ G_{n,i}^\ell$ the submanifold of
all $i$-dimensional subspaces $\xi$ with the property that all
canonical angles between $\xi$ and $\rl$ are equal.
\end{definition}

\begin{remark}\label{rm26} It is known, that the  Radon transform is overdetermined if the
dimension of the target space is greater than the dimension of the
source space. If $f$ is $\kl$-invariant and $i\le n-\ell$, then, by
Corollary \ref {cor32} and Remark \ref{Lzd}, the overdeterminicity
can be  eliminated if we restrict $(R_if)(\xi)$ to $\xi \in
G_{n,i}^\ell$.  Here one should mention the general method of the
kappa-operator, which allows us to reduce overdeterminicity by
invoking the relevant
 permissible complexes of subspaces; see, e.g., \cite {GGR} and
 references therein. The advantage of our {\it method of canonical angles}, which is applicable to
 the  particular case of $K_\ell$ -invariant functions, is the
 following. If  $i>\ell$, then to recover $f$ from $ R_if$, it
suffices to invert a simple Abel  integral (\ref {elli}). If $i\le
\ell$, then $f$  expresses through  $ R_if$ without any
integro-differential operations.
\end{remark}

\section{$K_\ell$ -Symmetric Bodies and Comparison of Volumes}

\subsection{Preliminaries}

An  origin-symmetric  (o.s.) star body $B$ in $\bbr^n$, $ n \ge 2$,
is a compact set with non-empty interior such that $tB \subset B \;
\forall t\in [0,1]$, $B=-B$, and  the {\it radial function} $ \rho_B
(\theta) = \sup \{ \lambda \ge 0: \, \lambda \theta \in B \}$ is
continuous on $S^{n-1}$.  The {\it Minkowski functional}  of $B$ is
defined by $ ||x||_B =\min \{a \ge 0 \, : \, x \in aB\}$, so that
$||\theta||_B=\rho_B^{-1}(\theta)$. An o.s.  star body $B$ is called
infinitely smooth if $\rho_B(\theta)\in C^\infty_{even}(S^{n-1})$.

If $\xi\in G_{n,i}, 1<i<n$, then \be\label{vol} \vol_i(B\cap \xi) =
i^{-1}\intl_{S^{n-1}\cap \xi} \rho_B^i (\theta) \,d_\xi
\theta=i^{-1}(R_i \rho_B^i )(\xi).\ee Similarly,
 $\vol_n (B) =
n^{-1}\int_{S^{n-1}} \rho_B^n (\theta) \,d\theta$.

\noindent{\bf Problem.} {\it Let $i$ be a  fixed integer, $1 \le i
\le n-1$. We wonder, for which o.s. star bodies $A$ and $B$  in
$\bbr^n$ the inequality  \be \label{bp1}\vol_i(A \cap \xi) \le
\vol_i(B \cap \xi) \quad \forall \xi \in G_{n,i} \ee implies
\be\label {bp2}\vol_n(A) \le \vol_n(B).\ee}

For $i=1$ an affirmative answer is obvious. Unlike the question in
Introduction, here we do not assume {\it a priori} that A and B are
convex. The reason is that  (\ref{bp1}) $\to$ (\ref{bp2}) may be
valid without any convexity assumption (see Theorem \ref{teor111}
(a)) and we want to understand how the convexity comes into play.

Below we study this problem when the body $A$ with smaller sections
is symmetric with respect to some mutually orthogonal subspaces,
say, $p$ and $p^\perp$, of dimensions $\ell$ and $n-\ell$,
respectively. We fix the coordinate system so that $
p=\rl=\overset{n}{\underset{j=n-\ell +1}{\oplus}}\,\bbr e_j$ and
$p^\perp =\rnl=\overset{n-\ell}{\underset{j=1}{\oplus}}\,\bbr e_j$.
Then $K_\ell A=A$, where $K_\ell$  is the group (\ref{kael}). An
o.s. star body with this property  is said to be $K_\ell$
-symmetric.

 By Lemma \ref {teo1}, the radial function $\rho_A
(\theta)$ of a $K_\ell$-symmetric   body $A$ is completely
determined by the angle $\om$ between the unit vector $\theta$ and
the subspace $\rl$. Hence,  we can set \be\label{rA} \rho_A
(\theta)= \tilde \rho_A (t), \qquad t=\cos^2 \om=\theta'{\rm
P}_{\rl}\theta.\ee By Theorem \ref{teo31}, the Radon transform $(R_i
f)(\xi), \; \xi \in  G_{n,i}$, of every $K_\ell$ -invariant function
$f$ is actually a function of the canonical angles between $\xi \in
G_{n,i}$ and $\rl$. Restricting $(R_i f)(\xi)$ to $\xi \in
G_{n,i}^\ell$ (see Definition \ref {defi}), we can remove
overdeterminicity of $R_i f$. As we shall see below, the lower
dimensional Busemann-Petty problem inherits this overdeterminicity,
and the latter can be eliminated in the same way by considering
sections by subspaces $\xi \in G_{n,i}^\ell$ only.

We will need the following auxiliary lemmas.
\begin{lemma}\label{50} The group $K_\ell$ preserves  canonical
angles between $\xi \in G_{n,i}$ and $\rl$.
\end{lemma}
\begin{proof} The proof relies on (\ref{forg}). Let first $\ell <i, \; \xi= \{\tau\}, \; \tau \in
V_{n,i}$. It suffices to check that for every $\gam \in K_\ell\,$,
matrices $r=\sig'\tau \tau'\sig$ and $r_\gam=\sig'\gam\tau
\tau'\gam'\sig$ have the same eigenvalues. Let $\gam=
\left[\begin{array} {cc} \a & 0
\\ 0 & \b \end{array} \right]$, where $\a  \!\in  \!O(n-\ell)$, $ \b  \!\in \!
O(\ell)$. Multiplying matrices, we have $ \gam'\sig=\sig\b'$. Hence,
$r_\gam=\b\sig'\tau \tau'\sig\b'=\b r\b'$. Since $\b r\b'$ and $r$
have the same eigenvalues,  we are done.

If $\ell \ge i$, we  compare  eigenvalues of matrices
$s=\tau'\sig\sig'\tau$ and $s_\gam=\tau'\gam'\sig\sig'\gam\tau$.
These matrices coincide, because, as we have already seen, $
\gam'\sig=\sig\b'$, and therefore,
$s_\gam=\tau'\sig\b'\b\sig'\tau=\tau'\sig\sig'\tau=s$.
\end{proof}
\begin{definition} {\rm($\kl$-symmetrization)} Given an o.s. star body $B$ in
$\bbr^n$, we define the associated $K_\ell$-symmetric body $B_0$ by
\be\label{ave} \rho_{B_0}(\theta)=\Big (\,\intl_{K_\ell} \rho_{B}^i
(\gam\theta)\, d\gam\Big )^{1/i}.\ee
\end{definition}
\begin{lemma}\label{510} $\vol_n (B_0)\le \vol_n (B)$.
\end{lemma}
\begin{proof}
By the generalized Minkowski inequality, \bea \vol_n^{i/n}
(B_0)&=&\Big (\frac {1}{n}\intl_{S^{n-1}}\Big [\intl_{K_\ell}
\rho_{B}^i (\gam\theta)\, d\gam\Big ]^{n/i}d\theta \Big
)^{i/n}\nonumber \\&\le&\intl_{K_\ell}\Big [\frac
{1}{n}\intl_{S^{n-1}}\rho_{B}^n (\gam\theta)\,d\theta \Big ]^{i/n}
d\gam=\vol_n^{i/n} (B),\nonumber \eea and the result follows.
\end{proof}
\begin{lemma}\label{cor45} Let $A$ and $B$ be o.s. star bodies  in
$\bbr^n$, $1\le \ell\le n-1$. If $A$ is $K_\ell$-symmetric, and
\be\label{ineq5} \vol_i(A \cap \xi) \le \vol_i(B \cap \xi)\quad
\forall \xi \in G_{n,i}^\ell,\ee then $ \vol_i(A \cap \xi) \le
\vol_i(B_0 \cap \xi)$ for all $\xi \in G_{n,i}^\ell$.
\end{lemma}
\begin{proof} Fix $\xi \in
G_{n,i}^\ell$. By Lemma \ref{50}, $\gam\xi \in G_{n,i}^\ell$ for
every $\gam \in K_\ell$. Owing to (\ref{ineq5}), $\vol_i(A \cap
\gam\xi) \le \vol_i(B \cap \gam\xi)$ or  $(R_i \rho_A^i
)(\gam\xi)\le(R_i \rho_B^i )(\gam\xi)$ for all $ \gam \in K_\ell$.
Integrating this inequality in $\gam$ and taking into account that
$R_i$ commutes with orthogonal transformations, we obtain
\be\label{55} (R_i \rho_A^i )(\xi)\le R_i \Big [ \intl_{K_\ell}
\rho_B^i (\gam\theta)\, d\gam \Big ](\xi)=(R_i \rho_{B_0}^i )(\xi).
\ee This implies  $ \vol_i(A \cap \xi) \le \vol_i(B_0 \cap \xi)$.
\end{proof}

\subsection{The case $  i\le \min(\ell, n-\ell)$} The following proposition represents
part (a) of  Theorem \ref{teor111}.
\begin{proposition} \label{pro1}Let $\;1\le i,\ell\le n-1$; $\; i\le \min(\ell, n-\ell)$.
If $\;A$ is a $K_\ell$-symmetric   body  in $\bbr^n$, then the
implication \be\label{bp3} \vol_i(A \cap \xi) \le \vol_i(B \cap \xi)
\;\forall \xi \in G_{n,i}^\ell \Longrightarrow \vol_n(A) \le
\vol_n(B)\ee is true for every o.s. star body $B$.
\end{proposition}
\begin{proof} By Lemma \ref{cor45}, for all $\xi \in
G_{n,i}^\ell$ we have  $$ \vol_i(A \cap \xi) \le \vol_i(B_0 \cap
\xi)\quad \text{\rm or}\quad (R_i \rho_A^i)(\xi)\le (R_i
\rho_{B_0}^i)(\xi).$$ Hence, by (\ref{gega}) and (\ref{rA}), $\tilde
\rho_A^i (\lam)\le \tilde \rho_{B_0}^i (\lam)$ for all $\lam \in
(0,1)$ (see Remark \ref{Lzd}), and therefore, $\rho_A (\theta)\le
\rho_{B_0} (\theta)$ for all $\theta \in S^{n-1}$. By Lemma \ref
{510}, it follows that $\vol_n(A)\le \vol_n (B_0)\le \vol_n (B)$.
\end{proof}

\subsection{The case $\ell<i\le n-\ell$}

We will need  some sort of duality which is a one-dimensional analog
of (\ref{dual})  and serves as a substitute for  the
 Lutwak's connection \cite {Lu} between the Busemann-Petty
problem and intersection bodies. According to (\ref{elli}), the
Radon transform $(R_i \rho_A^i)(\xi)$, restricted to $\xi \in
G_{n,i}^\ell$, is represented by the Abel type integral
\be\label{rhs}(R_i
\rho_A^i)(\xi)\!=\!\frac{c_1}{\lam^{i/2-1}}\intl_0^\lam \!t^{\ell/2
-1}(\lam\! -\!t)^{(i-\ell)/2 -1}\tilde\rho_A^i(t)\,
 dt, \ee
 $$ c_1=\sig_{i-\ell -1}\,\sig_{\ell
-1}/2,$$ where $\lam^{1/2}\in (0,1)$ is  the cosine of the
  canonical angles between $\xi$ and $\rl$ (we remind that these angles  are equal when $\xi \in G_{n,i}^\ell$
and (\ref{rhs}) is available for all $\lam \in (0,1)$; see Remark
\ref{Lzd}). Denote the right-hand side of (\ref{rhs}) by
$(I_+\tilde\rho_A^i)(\lam)$ and define the dual integral operator
\be\label{ala} (I_-\psi)(t)=c_1\, t^{\ell/2 -1}\intl_t^1 (\lam
-t)^{(i-\ell)/2 -1} \psi (\lam)\frac{d\lam}{\lam^{i/2-1}},\ee so
that \be\label{dual1} \intl_0^1 (I_+\tilde\rho_A^i)(\lam)\,\psi
(\lam)\,d\lam=\intl_0^1 \tilde\rho_A^i(t)\,(I_-\psi)(t)\, dt.\ee
Expression (\ref {ala}) resembles the classical Riemann-Liouville
integral \be (I_-^\a g)(t)=\frac{1}{\Gam (\a)}\intl_t^1 g(\lam)(\lam
-t)^{\a -1}\, d\lam, \qquad \a>0.\ee
\begin{lemma}\label{lem} Let $\;1\le \ell<i\le n-\ell$ and suppose that $A$ is a
$K_\ell$ -symmetric   body  in $\bbr^n$. If there is a non-negative
function $g$ on $(0,1)$, which is integrable on every interval
$(\del,1)$, $0<\del <1$, and such that \be\label{cond}
(1-t)^{(n-\ell)/2-1}\tilde \rho_A^{n-i} (t)=(I_-^\a g)(t), \qquad
\a=\frac{i-\ell}{2}, \quad t \in (0,1),\ee
 then the implication \be\label{bp4} \vol_i(A \cap \xi) \le
\vol_i(B \cap \xi) \;\forall \xi \in G_{n,i}^\ell \Longrightarrow
\vol_n(A) \le \vol_n(B)\ee holds for every o.s. star body $B$.
\end{lemma}
\begin{proof} By (\ref{gn1}),
$$\vol_n(A)=\frac{1}{n}\,\intl_{S^{n-1}}\rho_A^n (\theta)\,
d\theta=c_2 \intl_0^{1}\tilde \rho_A^n (t) t^{\ell/2-1}
(1-t)^{(n-\ell)/2-1}\, dt,$$ $$c_2=\sig_{\ell -1}\sig_{n-\ell
-1}/2n.$$ Hence, owing to  (\ref{ala}), (\ref{dual1}), and
(\ref{cond}), \bea\vol_n(A)&=&c_2 \intl_0^{1}\tilde \rho_A^i (t)
t^{\ell/2-1}(I_-^\a g)(t)\, dt=\frac{c_2}{c_1}\intl_0^{1}\tilde
\rho_A^i (t)(I_-\psi)(t)\, dt\nonumber \\
&=&\frac{c_2}{c_1}\intl_0^{1}(I_+\tilde\rho_A^i)(\lam)\,\psi
(\lam)\,d\lam, \qquad \psi (\lam)\!=\! \frac{\lam^{i/2-1}
g(\lam)}{\Gam (\a)}\ge 0.\eea If $\vol_i(A \cap \xi) \le \vol_i(B
\cap \xi) \;\forall \xi \in G_{n,i}^\ell$, then, by (\ref{rhs}) and
(\ref{55}),
$$
(I_+\tilde\rho_A^i)(\lam)=(R_i \rho_A^i)(\xi)\le (R_i \rho_{B_0}^i
)(\xi)=(I_+\tilde\rho_{B_0}^i)(\lam),
$$
and therefore, \bea \vol_n(A)&\le&
\frac{c_2}{c_1}\intl_0^{1}(I_+\tilde\rho_{B_0}^i)(\lam)\,\psi
(\lam)\,d\lam\nonumber \\
&=&c_2\intl_0^1 \tilde\rho_{B_0}^i (t)\tilde \rho_A^{n-i} (t)
t^{\ell/2-1} (1-t)^{(n-\ell)/2-1}\, dt\nonumber\\
&=&\frac{1}{n}\,\intl_{S^{n-1}}\rho_{B_0}^i (\theta)\rho_A^{n
-i}(\theta)\, d\theta.\nonumber \eea Now H\"older's inequality
yields $\vol_n(A)\le \vol_n(B_0)$, and the result follows by Lemma
\ref{510}.
\end{proof}

Up to now,  the  $K_\ell$-symmetric body $A$ with smaller sections
was  arbitrary. To handle the case $i>\ell$, we additionally assume
that $A$ is convex. The following  lemma enables us to reduce
consideration to smooth bodies.

\begin{lemma}\label{arbi} Let $A$ and $B$ be  o.s. star bodies in $\bbr^n$.
 If the implication
\be\label{bp4} \vol_i(A \cap \xi) \le \vol_i(B \cap \xi) \;\forall
\xi \in G_{n,i}^\ell \Longrightarrow \vol_n(A) \le \vol_n(B)\ee is
true for every infinitely smooth $K_\ell$-symmetric convex body $A$,
then it is true when $A$ is an arbitrary $K_\ell$-symmetric convex
body.
\end{lemma}
\begin{proof} Given a $K_\ell$-symmetric convex body
$A$, let $A^*=\{x: |x\cdot y| \le 1\; \forall y\in A \}$ be the
polar body of  $A$ with the support function $h_{A^*} (x)=\max \{
x\cdot y : y \in A^*\}$. Since $h_{A^*} (\cdot)$ coincides with
Minkowski's functional $ ||\cdot||_A$, then $h_{A^*} (\cdot)$ is
$\kl$-invariant, and therefore, $A^*$ is $K_\ell$-symmetric. It is
known \cite [pp. 158-161] {Schn}, that any o.s. convex body in $\rn$
can be approximated   by infinitely smooth convex bodies with
positive curvature and the approximating operator commutes with
rigid motions. Hence, there is a sequence $\{A_j^*\}$ of infinitely
smooth $K_\ell$-symmetric convex bodies with positive curvature such
that $h_{A^*_j} (\theta)$ converges to $h_{A^*}(\theta)$ uniformly
on $S^{n-1}$. The latter means, that for the relevant sequence of
infinitely smooth $K_\ell$-symmetric convex bodies $A_j=(A_j^*)^*$,
$$\lim_{j \to \infty }\,\max_{\theta \in S^{n-1}} |\, ||\theta||_{A_j}-
||\theta||_A |=0.$$
 This implies convergence in the radial metric, i.e.,
\be\label{conv1}\lim_{j \to \infty }\,\max_{\theta \in S^{n-1}}
|\rho_{A_j}(\theta) - \rho_{A}(\theta)|=0.\ee

Let us show that the sequence $\{A_j\}$ in (\ref{conv1}) can be
modified so that $A_j \subset A$. The idea of this argument was
borrowed from \cite {RZ}. Without loss of generality, assume that
$\rho_A(\theta)\ge1$. Choose $A_j$ so that
$$
\vert \rho_{A_j}(\theta)-\rho_A(\theta)\vert < \frac1{j+1} \quad
\forall \theta \in S^{n-1}
$$
and set  $A_j'=\frac j{j+1} A_j$. Then, obviously,
$\rho_{A_j'}(\theta) \to \rho_A(\theta)$ uniformly on $S^{n-1}$ as
$j \to \infty$, and
$$
\rho_{A_j'}=\frac j{j+1} \rho_{A_j} <\frac j{j+1} \big
(\rho_A+\frac1{j+1}\big )\le \rho_A.
$$
Hence, $A_j'\subset A$.

Now suppose that  $\vol_i(A \cap \xi) \le \vol_i(B \cap \xi)
\;\forall \xi \in G_{n,i}^\ell$. Then this is  true when $A$ is
replaced by $A_j'$, and, by the assumption of the lemma,
$\vol_n(A_j') \le \vol_n(B)$. Passing to the limit as $j \to
\infty$, we obtain $\vol_n(A) \le \vol_n(B)$.
\end{proof}

The next proposition gives part (b) of  Theorem \ref{teor111}.
\begin{proposition}\label{pro2} Let $A$ be a
$K_\ell$-symmetric  convex body  in $\bbr^n$, and let $\;2\le i\le
n-\ell$. If

$$ i= \ell+1 \qquad  \text{\rm $($in this case $\ell \le (n-1)/2)$}$$
or
$$ i=\ell+2 \qquad  \text{\rm  $($in
this case $\ell \le (n-2)/2)$},$$
 then the implication \be\label{bp4}
\vol_i(A \cap \xi) \le \vol_i(B \cap \xi) \;\forall \xi \in
G_{n,i}^\ell \Longrightarrow \vol_n(A) \le \vol_n(B)\ee holds for
every o.s. star body $B$.
\end{proposition}
\begin{proof}  By Lemma \ref{arbi}, we can assume $\rho_A
\in C^\infty (S^{n-1})$. If $i= \ell+2$, then (\ref{cond}) becomes $
(1-t)^{(n-i)/2}\tilde \rho_A^{n-i} (t)=\int_t^1 g(\lam)\, d\lam$,
which implies
$$g(t)=-\frac{d}{dt}[(1-t)^{(n-i)/2}\tilde \rho_A^{n-i} (t)]\in L^1
(0,1).$$ To check that $g$ is nonnegative, we set $t=1-s$, $
r(s)=s^{1/2}\tilde \rho_A (1-s)$, $s=\sin^2 \om$, and get
$$ g(1-s)=\frac{d}{ds}[r^{n-i}(s)]=(n-i)r^{n-i-1}(s)r'(s).
$$
If $\theta=u\,\sin \om +v\,\cos\, \om\in S^{n-1}$, $u\in
S^{n-\ell-1} \subset \bbr^{n-\ell}, \; v \in S^{\ell -1}\subset\rl$,
and $P_{u,v}$ is a 2-plane spanned by  $u$ and $v$, then $A\cap
P_{u,v}$ is a convex domain, which is symmetric with respect to the
$u$ and $v$ axes. Since $s=\sin^2 \om$, then  $r(s)=s^{1/2}\tilde
\rho_A(1-s)$ is non-decreasing, and therefore, $r'(s)\ge 0$. This
gives $g(1-s)\ge 0$, $s \in (0,1)$, or, equivalently, $g(t) \ge 0$
for all $0<t<1$. Now the result follows by Lemma \ref{lem}.

Let $i= \ell+1$. We set $\varkappa_A (t)=(1-t)^{(n-i-1)/2}\tilde
\rho_A^{n-i} (t)$ and reconstruct $g(t)$  from (\ref{cond}) using
fractional differentiation as follows: \bea g(t)&=&-\frac{1}{\sqrt
{\pi}}\,\frac{d}{dt}\intl_t^1 (s -t)^{ -1/2}\varkappa_A (s)\, ds
\quad \text{\rm (set $p=1-t, \;
q=1-s$)}\nonumber\\
&=&\frac{1}{\sqrt {\pi}}\,\frac{d}{dp}\intl_0^p (p-q)^{
-1/2}\varkappa_A (1-q)\, dq\nonumber\\
&=&\frac{1}{\sqrt {\pi}}\,\frac{d}{dp}\Big [p^{1/2}\intl_0^1
(1-\eta)^{ -1/2}\varkappa_A (1-p\eta)\, d\eta\Big ].\nonumber\eea
This gives
 \bea g(t)
&=& \frac{1}{\sqrt {\pi}}\,\frac{d}{dp}\intl_0^1
[(p\eta)^{1/2}\tilde
\rho_A (1-p\eta)]^{n-i}\frac{d\eta}{\sqrt {\eta (1-\eta)}}\,d\eta\nonumber\\
&=& \frac{1}{\sqrt {\pi}}\,\frac{d}{dp}\intl_0^1
\frac{r^{n-i}(p\eta)}{\sqrt {\eta (1-\eta)}}\,d\eta, \qquad
r(s)=s^{1/2}\tilde \rho_A (1-s).\nonumber\eea The last integral is a
non-decreasing function of $p$, and therefore, the derivative of it
is non-negative. Hence, $g(t) \ge 0$ for all $0<t<1$ and, by Lemma
\ref{lem}, we are done.
\end{proof}

\section{The negative result}

The proof of the negative result in  Theorem \ref{teor112} relies on
Koldobsky's  generalizations of the Lutwak's  concept of
intersection body (see \cite{K4}, \cite{Lu}) and properties of the
generalized cosine transforms \cite{R3}.

We remind basic facts. The  generalized cosine transform of a
function $f$ on $S^{n-1}$ is defined by \be\label{af7} (M^\a f)(u)=
\gam_n(\a)\, \int_{S^{n-1}} f(\theta) |\theta \cdot u|^{\a-1}
\,d\theta, \ee \be\label{beren7} \gam_n(\a)\!=\!{
\sig_{n-1}\,\Gamma\big( (1\!-\!\a)/2\big)\over 2\pi^{(n-1)/2} \Gamma
(\a/2)}, \qquad Re \, \a \!>\!0, \quad \a \!\neq \!1,3,5, \ldots
.\ee The integral (\ref{af7}) is absolutely convergent for any $f\in
L^1(S^{n-1})$. If $f$ is infinitely differentiable, then $M^\a f $
extends as meromorphic function of $\a$  with the poles $\a=1,3,5,
\ldots$. The following statement is a consequence of the relevant
spherical harmonic decomposition.
\begin{lemma}\label{l17} \cite{R5} Let $\a, \b \in \bbc; \; \a, \b \neq
1,3,5, \ldots \,$. If $\a+\b=2-n$ and $f\in
C^\infty_{even}(S^{n-1})$, then \be\label{st7}M^\a M^{\b}f=f.\ee If
$\, \a, 2-n-\a \neq 1,3,5, \ldots $, then $M^\a$ is an automorphism
of the space $C^\infty_{even}(S^{n-1})$ endowed with the standard
topology.
\end{lemma}

We will also need the following statement, which is a particular
case of Lemma 3.5 from \cite{R3}.
\begin{lemma}\label{l27} Let $f \in C^\infty_{even}(S^{n-1})$. Then
\be \label{con7} (R_i M^{1-i} f)(\xi)=c\,
 (R_{n-i}f)(\xi^\perp), \qquad \xi \in G_{n,i},\ee where $c=c(n,i)$ is a
positive constant.
\end{lemma}
\begin{definition}\label{def47}  \cite{K4} An o.s. star body $K$ in $\bbr^n$ is  a
 $k$-intersection body if there is  a
non-negative finite Borel measure $\mu$ on $S^{n-1}$, so that for
every Schwartz function $\phi$,
\[\int_{\bbr^n}||x||_K^{-k} \phi (x)\, dx=\int_{S^{n-1}}\Big
[\int_0^\infty t^{k-1}\hat \phi (t\theta) \,dt \Big ] \,d\mu
(\theta),\] where $\hat \phi$ is the Fourier transform of $\phi$. We
denote by $\I^n_k$ the class of all $k$-intersection bodies in
$\bbr^n$.
\end{definition}

The following equivalent  definition is a particular case of the
more general Definition 5.4 from \cite{R3}.
\begin{definition}\label{def48}   An o.s. star body $K$ in $\bbr^n$ is  a
 $k$-intersection body if there is  a
non-negative finite Borel measure $\mu$ on $S^{n-1}$, so that
$\rho_K^k=M^{1-k}\mu$, i.e., $(\rho_K^k, \vp)=(\mu, M^{1-k}\vp)$ for
any $\vp\in C^\infty_{even}(S^{n-1})$.
\end{definition}

The next proposition plays a key role in the proof of the negative
result in this section.
\begin{lemma}\label{finn7}
Let $B$ be an infinitely smooth $K_\ell$ -symmetric  convex body
with positive curvature. If $B \notin \I^n_{n-i}$, then there is an
infinitely smooth $K_\ell$-symmetric  convex body $A$ such that $
\vol_i(A \cap \xi) \le \vol_i(B \cap \xi)\quad \forall \xi \in
G_{n,i}$, but $\vol_n(A)
> \vol_n(B)$.
\end{lemma}
\begin{proof} By Definition \ref{def48} with $k=n-i$,  there is a
function  $\vp$ in $C^\infty_{even}(S^{n-1})$, which is negative on
some open
 set $\Om
 \subset S^{n-1}$ and such that $ \rho_B^{n-i}=M^{1+i-n} \vp$.
Since $B$ is $K_\ell$ -symmetric and operators $M^{\a}$ commute with
orthogonal transformations, then $\vp$ is $K_\ell$ -invariant and
$\vp<0$ on the whole orbit $\Om_\ell=K_\ell \Om$. Choose a function
$h\in C^\infty_{even}(S^{n-1})$ so that $h \not\equiv 0$,
$h(\theta)\ge 0$ if $\theta \in \Om_\ell$ and $h(\theta)\equiv 0$
otherwise. Without loss of generality, we can assume $h$ to be
$K_\ell$ -invariant (otherwise, we can replace it by $\tilde
h(\theta)=\int_{K_\ell}h(\gam\theta) d\gam$). Define an
origin-symmetric smooth  body $A$ by $ \rho_A^{i}=\rho_B^{i}-\e
M^{1-i} h$, $\e>0$. Clearly, $A$ is $K_\ell$ -symmetric. If $\e$ is
small enough, then $A$ is convex. This conclusion is a consequence
of Oliker's formula \cite{Ol}, according to which the Gaussian
curvature of an o.s. star body expresses through the first and
second derivatives of the radial function. Since by (\ref{con7}),
$(R_i M^{1-i} h)(\xi)=c\,
 (R_{n-i}h)(\xi^\perp)\ge 0$, then $R_i\rho_A^{i}\le R_i\rho_B^{i}$. This
 gives $\vol_i(A \cap \xi) \le \vol_i(B \cap \xi)\quad \forall \xi \in
G_{n,i}$.  On the other hand, by (\ref{st7}),
$$ ( \rho_B^{n-i},\rho_B^{i} -\rho_A^{i})=\e(M^{1+i-n} \vp, M^{1-i}
h)=\e(\vp, h)<0,$$ or $( \rho_B^{n-i},\rho_B^{i})<(
\rho_B^{n-i},\rho_A^{i})$. By H\"older's inequality, this implies
 $ \vol_n(B)<\vol_n(A)$.
\end{proof}

Consider the $(q,\ell)$-ball $ B^n_{q,\ell}=\{x =(x', x''): |x'|^q
+|x''|^q \le 1\}$, where $x' \in
\rnl=\overset{n-\ell}{\underset{j=1}{\oplus}}\,\bbr e_j,    \quad
x''\in \rl=\overset{n}{\underset{j=n-\ell +1}{\oplus}}\,\bbr e_j$.
\begin{lemma}\label{thm217} $($cf. \cite [Theorem 4.21]{K3}$)$ Let
 \be\label{ballm} B^{m+1}_{q,1}=\{(x',y): |x'|^q +|y|^q \le 1,
\; x'\in \bbr^m, y \in \bbr\}.\ee If $q>2$ and $m\ge k+3$, then
$B^{m+1}_{q,1}\notin \I^{m+1}_k$.
\end{lemma}
\begin{lemma}\label{nakz}
If $q>2$ and $\ell +2<i\le n-1$, then
 $B^n_{q,\ell}\notin \I^n_{n-i}$.
\end{lemma}
\begin{proof}  Suppose the contrary and
consider the section of $B^n_{q,\ell}$ by the $(n-\ell
+1)$-dimensional  plane $\eta=\bbr e_n\oplus\rnl$. By Proposition
3.17 from \cite{Mi2} (see also more general Theorem 5.12 in
\cite{R3}) $B^n_{q,\ell} \cap \eta \in \I^{n-\ell +1}_{n-i}$ in
$\eta$, but this contradicts Lemma \ref{thm217}, according to which
(set $m=n-\ell$) $B^n_{q,\ell} \cap \eta$ is not an
$(n-i)$-intersection body when $i>\ell +2$.
\end{proof}

For $q=4$, Lemmas \ref{finn7} and \ref{nakz} imply the following
negative result.
\begin{proposition} If  $\ell +2<i\le n-1$, then
 there is an infinitely smooth $K_\ell$-symmetric
convex body $A$ such that $$ \vol_i(A \cap \xi) \le
\vol_i(B^n_{4,\ell} \cap \xi)\quad \forall \xi \in G_{n,i}, \quad
\text { but}\quad \vol_n(A)
> \vol_n(B^n_{4,\ell}).$$
\end{proposition}
This is just  Theorem \ref{teor112}.

\end{document}